\documentclass[10pt]{amsart}
\usepackage{mathrsfs}
\usepackage{amsfonts} 
\usepackage{graphicx,latexsym,bm,amsmath,amssymb,verbatim,multicol,lscape}
\theoremstyle{definition}

\theoremstyle{remark}

\numberwithin{equation}{section}

\begin{document}
\title[Necessary conditions for RDPS to be permutational]
{Necessary conditions for reversed Dickson polynomials of the second kind to be permutational}%
\author[Shaofang Hong and Xiaoer Qin]{Shaofang Hong and Xiaoer Qin\\
Mathematical College, Sichuan University, Chengdu 610064, P.R. China}
\thanks{The research was supported partially by National Science Foundation
of China Grant \#11371260. Emails: sfhong@scu.edu.cn, s-f.hong@tom.com,
hongsf02@yahoo.com (S. Hong). qincn328@sina.com (X. Qin)}
\date{\today}
\begin{abstract}
In this paper, we present several necessary conditions for
the reversed Dickson polynomial $E_{n}(1, x)$ of the second kind
to be a permutation of $\mathbb{F}_{q}$. In particular, we give
explicit evaluation of the sum $\sum_{a\in \mathbb{F}_{q}}E_{n}(1,a)$.
\end{abstract}

\maketitle

\section{\bf Introduction}

Let $p$ be a prime and $\mathbb{F}_{q}$ be a finite field of $q=p^{e}$ elements,
where $e$ is a positive integer. Associated to any integer $n\geq 0$ and a parameter
$a\in\mathbb{F}_{q}$, the $n$-th {\it Dickson polynomials of the first kind and of
the second kind}, denoted by $D_{n}(x,a)$ and $E_{n}(x,a)$, are defined by
$$D_{n}(x,a):=\sum_{i=0}^{[\frac{n}{2}]}\frac{n}{n-i}
\binom{n-i}{i}(-a)^{i}x^{n-2i}$$
and
$$E_{n}(x, a):=\sum_{i=0}^{[\frac{n}{2}]}\binom{n-i}{i}(-a)^{i}x^{n-2i},$$ respectively.
Recently, Wang and Yucas \cite {[WY]} further defined the {\it $n$-th Dickson polynomial of
the $(k+1)$-th kind} $D_{n,k}(x,a)\in \mathbb{F}_{q}[x]$ by
$$D_{n,k}(x,a):=\sum_{i=0}^{[\frac{n}{2}]}\frac{n-ki}{n-i}
\binom{n-i}{i}(-a)^{i}x^{n-2i}.$$
On the other hand, Hou, Mullen, Sellers and Yucas \cite{[HMSY]} introduced
the definition of {\it the reversed Dickson polynomial of the first kind},
denoted by $D_{n}(a,x)$, as follows
$$D_{n}(a,x):=\sum_{i=0}^{[\frac{n}{2}]}\frac{n}{n-i}
\binom{n-i}{i}(-x)^{i}a^{n-2i}.$$ By extending the definition of reversed
 Dickson polynomials, Wang and Yucas \cite {[WY]} got the definition
of {\it the $n$-th reversed Dickson
polynomial of the $(k+1)$-th kind} $D_{n,k}(a,x)\in \mathbb{F}_{q}[x]$,
which is defined by 

$$D_{n,k}(a,x):=\sum_{i=0}^{[\frac{n}{2}]}\frac{n-ki}{n-i}
\binom{n-i}{i}(-x)^{i}a^{n-2i}.$$

The permutation behavior of
Dickson polynomials $D_{n}(x,a)$ over finite fields are well known:
$D_{n}(x,0)=x^{n}$ is a permutation polynomial of $\mathbb{F}_{q}$
if and only if $(n,q-1)=1$, and if $a\neq0$, then $D_{n}(x,a)$
induces a permutation of $\mathbb{F}_{q}$ if and only if
$(n,q^{2}-1)=1$ (see \cite{[LN]}, Theorem 7.16). Meanwhile, there are many
results on permutation properties of Dickson polynomial $E_{n}(x,a)$
of the second kind, the readers can be referred to \cite{[C]}. In
\cite{[WY]}, Wang and Yucas studied the permutational behavior of
Dickson polynomials of the third kind $D_{n,2}(x,1)$. They obtained
some necessary conditions for $D_{n,2}(x,1)$ to be a permutation
polynomial of $\mathbb{F}_{q}$.

Hou, Mullen, Sellers and Yucas \cite{[HMSY]} studied the permutation
properties of reversed Dickson polynomial $D_{n}(a,x)$ of the first
kind. In fact, they showed that $D_{n}(a,x)$ is closely related to almost
perfect nonlinear (APN) functions, and got several families of permutation
polynomials from reversed Dickson polynomials of the first kind.
In \cite{[HL]}, Hou and Ly found several necessary conditions for
reversed Dickson polynomials $D_{n}(1,x)$ of the first kind
to be a permutation polynomial.

In this paper, we mainly investigate reversed Dickson polynomial of the
second kind. We denote by $E_{n}(a,x)\in \mathbb{F}_{q}[x]$ the reversed
Dickson polynomial of the second kind, which is defined by
$$E_{n}(a,x):=\sum_{i=0}^{[\frac{n}{2}]}\binom{n-i}{i}(-x)^{i}a^{n-2i}.\eqno(1.1)$$
For $a\neq 0$, we write $x=y(a-y)$ with an indeterminate $y\neq
\frac{a}{2}$.
Then $E_{n}(a,x)$ can be rewritten as
$$E_{n}(a,x)=\frac{y^{n+1}-(a-y)^{n+1}}{2y-a}.\eqno(1.2)$$
We will emphasize on the permutation behavior of reversed Dickson
polynomials $E_{n}(a,x)$ of the second kind over $\mathbb{F}_{q}$.
This paper is organized as follows. First in Section 2, we study the
properties of the reversed Dickson polynomial $E_{n}(a,x)$ of the second
kind. Consequently, in Section 3, by introducing the polynomial
$f_{m}(x)=\sum_{j=0}^{[\frac{m-1}{2}]}\binom{m}{2j+1}x^{j}$,
we prove several necessary conditions for
the reversed Dickson polynomial $E_{n}(1,x)$ of the second kind to
be a permutation polynomial of $\mathbb{F}_{q}$. It is well
known that a function $f:\mathbb{F}_{q}\to \mathbb{F}_{q}$
is a permutation polynomial of $\mathbb{F}_{q}$ if and only if
\[\sum_{a\in\mathbb{F}_{q}}f(a)^i=
 \begin{cases}\ \ 0,\ \ \ {\rm if} &0\leq i\leq q-2,\\
-1,\ \ \ {\rm if} & i=q-1.
\end{cases}\]
Thus we would like to know if the sum
$\sum_{a\in \mathbb{F}_{q}}E_{n}(1,a)^i$ is computable.
We are able to treat with this sum when $q$ is odd and $i=1$.
The final section is devoted to the computation of the sum
$\sum_{a\in \mathbb{F}_{q}}E_{n}(1,a)$.

\section{\bf Reversed Dickson polynomials of the second kind}

In this section, we mainly study properties of reversed
Dickson polynomials of the second kind. If $a=0$, then
\[E_{n}(0,x)=\begin{cases}\ \ \ 0,\ \  \ \ \ {\rm if} \  n\ {\rm is\ odd},\\
(-x)^{k}, \ \ {\rm if}\ n=2k, k\ {\rm is\ nonnegative\ integer}.
\end{cases}\]
Hence $E_{n}(0,x)$ is a PP (permutation polynomial) of $\mathbb{F}_{q}$ if
and only $n=2k$ with $(k,q-1)=1$. In what follows we assume that $a\in\mathbb{F}_{q}^{\ast}$.
By a trivial fact that $f(x)$ is a PP of $\mathbb{F}_{q}$ if and only if $cf(dx)$
is a PP of $\mathbb{F}_{q}$ for any given $c, d \in\mathbb{F}_{q}^*$, we can easily
deduce the following result.

\noindent{\bf Theorem 2.1.} {\it Let $a,b\in\mathbb{F}_{q}^{\ast}$.
Then $E_{n}(a,x)=\frac{a^{n}}{b^{n}}E_{n}\big(b,\frac{b^{2}}{a^{2}}x\big).$
Furthermore, $E_{n}(a,x)$ is a PP of $\mathbb{F}_{q}$ if and only if
$E_{n}(1, x)$ is a PP of $\mathbb{F}_{q}$.
}
\begin{proof}
First by the definition (1.1), we have 
\begin{eqnarray*}
\frac{a^{n}}{b^{n}}E_{n}\Big(b,\frac{b^{2}}{a^{2}}x\Big)&=&\frac{a^{n}}{b^{n}}
\sum_{i=0}^{[\frac{n}{2}]}\binom{n-i}{i}\Big(-\frac{b^{2}}{a^{2}}x\Big)^{i}b^{n-2i}\\
&=&\sum_{i=0}^{[\frac{n}{2}]}\binom{n-i}{i}(-x)^{i}b^{n-2i}\Big(\frac{b}{a}\Big)^{2i-n}\\
&=&E_{n}(a,x).
\end{eqnarray*}
So the first part is proved.

To show the second part, one notices that $E_{n}(a,x)=a^{n}E_{n}(1,\frac{x}{a^{2}})$.
Since $a\in\mathbb{F}_{q}^{\ast}$, one has
$a^n, \frac{1}{a^2}\in\mathbb{F}_{q}^{\ast}$. It follows that
$E_{n}(a,x)$ is a PP of $\mathbb{F}_{q}$ if and only $E_{n}(1,x)$
is a PP of $\mathbb{F}_{q}$. This concludes the proof of second part.
Hence Theorem 2.1 is proved.
\end{proof}

By Theorem 2.1, it is easy to see that to study the permutation
behavior of reversed Dickson polynomial $E_{n}(a,x)$ of the
second kind, one needs only to consider that of $E_{n}(1, x)$.
In the following, we list some basic facts about the
reversed Dickson polynomial $E_{n}(1,x)$ of the second kind.

\noindent{\bf Theorem 2.2.} {\it Let $p$ be an odd prime, $n$ and $r$
be positive integers. Each of the following is true:\\
{\rm (1)}. We have $E_{n}(1,x(1-x))=\frac{x^{n+1}-(1-x)^{n+1}}{2x-1}$ if
$x\ne \frac{1}{2}$ and $E_{n}(1,\frac{1}{4})=\frac{n+1}{2^{n}}$.\\
{\rm (2)}. If $\gcd(p,n)=1$, then $E_{np^r-1}(1,x)
=(E_{n-1}(1,x))^{p^r}(1-4x)^{\frac{p^r-1}{2}}$. \\
{\rm (3)}. If $n_{1}$ and $n_{2}$ are positive integers such that $n_{1}\equiv
n_{2}\pmod {q^2-1}$, then $E_{n_{1}}(1,x_0)=E_{n_{2}}(1,x_0)$
for any $x_0\in \mathbb{F}_{q}\setminus\{\frac{1}{4}\}$.}
\begin{proof}
{\rm (1).} Clearly, the first identity follows from (1.2).
To prove the second one, we notice that by (2.3) of \cite{[C]}
(see page 226 of \cite{[C]}), we have $E_{n}(2,1)=n+1$. But
Theorem 2.1 tells us that $E_{n}(2,1)=2^{n}E_{n}(1,\frac{1}{4})$.
Thus $E_{n}(1,\frac{1}{4})=\frac{n+1}{2^{n}}$ as required.

{\rm (2).} Writing $x=y(1-y)$ with $y\ne \frac{1}{2}$ being an
indeterminate gives us that $1-4x=(2y-1)^2$. So by part (1),
one derives that
\begin{eqnarray*}
E_{np^r-1}(1,x)&=&E_{np^r-1}(1,y(1-y))\\
&=&\frac{y^{np^r}-(1-y)^{np^r}}{2y-1}\\
&=&\Big(\frac{y^{n}-(1-y)^{n}}{2y-1}\Big)^{p^r}(2y-1)^{p^r-1}\\
&=&\big(E_{n-1}(1,y(1-y))\big)^{p^r}(2y-1)^{p^r-1}\\
&=&E_{n-1}(1,x)^{p^r}(1-4x)^{\frac{p^r-1}{2}}.
\end{eqnarray*}
Particularly, if $x=\frac{1}{4}$, then by part (1), we have
$$E_{np^r-1}(1,x)=E_{np^r-1}\Big(1,\frac{1}{4}\Big)=\frac{np^r}{2^{np^r-1}}
=0=(E_{n-1}(1,x))^{p^r}(1-4x)^{\frac{p^r-1}{2}}$$
as desired. Part (2) is proved.

{\rm (3).} For each $x_0\in\mathbb{F}_{q}\setminus\{\frac{1}{4}\}$, one may write
$y_0\in\mathbb{F}_{q^{2}}\setminus\{\frac{1}{2}\}$ such that $x_0=y_0(1-y_0)$. Thus
\begin{eqnarray*}
E_{n_{1}}(1,x_0)&=&E_{n_{1}}(1,y_0(1-y_0))\\
&=&\frac{y_0^{n_{1}+1}-(1-y_0)^{n_{1}+1}}{2y_0-1}\\
&=& \frac{y_0^{n_{2}+1}-(1-y_0)^{n_{2}+1}}{2y_0-1} \\
&=&E_{n_{2}}(1,x_0).
\end{eqnarray*}
This ends the proof of Theorem 2.2.
\end{proof}

\noindent{\bf Remark.} When $p=2$, we have
$$E_{n}(1,x(1-x))=x^{n+1}+(1-x)^{n+1}=D_{n+1}(1,x(1-x)).$$
In \cite{[HMSY]}, Hou et al. discussed some connections between
reversed Dickson PPs of $\mathbb{F}_{q}$ and APN functions of
$\mathbb{F}_{q}$, and obtained several families of reversed
Dickson PPs. Throughout the reminder of this article, unless
specified, $p$ is always assumed to be an odd prime.

By \cite{[WY]}, we know that $E_{n}(x,a)=xE_{n-1}(x,a)-aE_{n-2}(x,a)$
holds for any integer $n\geq2$. Regarding $E_{n}(1,x)$, we have the following result.

\noindent{\bf Proposition 2.1.} {\it Let $p$ be an odd prime
and $n\ge 2$ be an integer. Then
$E_{n}(1,x)=E_{n-1}(1,x)-xE_{n-2}(1, x).$}
\begin{proof}
First we consider the case $x\ne \frac{1}{4}$. For this case, one may let $x=y(1-y)$
with $y$ being an indeterminate and $y\ne \frac{1}{2}$. Then by Theorem 2.2 (1), we have
\begin{eqnarray*}
&&E_{n-1}(1,y(1-y))-y(1-y)E_{n-2}(1, y(1-y))\\
&=&\frac{y^n-(1-y)^n}{2y-1}-y(1-y)\frac{y^{n-1}-(1-y)^{n-1}}{2y-1}\\
&=&\frac{y^n-(1-y)^n}{2y-1}-\frac{y^{n}(1-y)-y(1-y)^{n}}{2y-1}\\
&=&\frac{y^{n+1}-(1-y)^{n+1}}{2y-1}=E_{n}(1,y(1-y)).
\end{eqnarray*}

For the case $x=\frac{1}{4}$, by Theorem 2.2 (1), we infer that
$$
E_{n-1}\Big(1, \frac{1}{4}\Big)-\frac{1}{4}E_{n-2}\Big(1, \frac{1}{4}\Big)
=\frac{n}{2^{n-1}}-\frac{n-1}{2^{n}}=\frac{n+1}{2^{n}}=E_{n}\Big(1,\frac{1}{4}\Big).
$$
Thus Proposition 2.2 is proved.
\end{proof}

Using this recursion, we can obtain the generating function of the
reversed Dickson polynomial $E_{n}(1,x)$ of the second kind as follows.

\noindent{\bf Proposition 2.2.} {\it The generating function of $E_{n}(1,x)$
is given by:}
$$\sum_{n=0}^{\infty}E_{n}(1,x)t^{n}=\frac{1}{1-t+xt^{2}}.$$

\begin{proof} By Proposition 2.1, we have

\begin{eqnarray*}
&&(1-t+xt^{2})\sum_{n=0}^{\infty}E_{n}(1,x)t^{n}\\
&=&\sum_{n=0}^{\infty}E_{n}(1,x)t^{n}-\sum_{n=0}^{\infty}E_{n}(1,x)t^{n+1}
+x\sum_{n=0}^{\infty}E_{n}(1,x)t^{n+2}\\
&=&1+t-t+\sum_{n=0}^{\infty}(E_{n+2}(1,x)-E_{n+1}(1,x)+xE_{n}(1,x))t^{n+2}=1.
\end{eqnarray*}
Thus the desired result follows immediately.
\end{proof}

In the following, by using the reversed Dickson polynomial $E_{n}(1,x)$
of the second kind, we obtain some PPs of $\mathbb{F}_{q}$.

\noindent{\bf Proposition 2.3.} {\it Let $p$ be an odd prime and
$k$ be a positive integer. Then we have}
$$E_{p^{k}-1}(1,x)=(1-4x)^{\frac{p^{k}-1}{2}}.$$
\begin{proof}
First putting $x=y(1-y)$ with an indeterminate $y\neq
\frac{1}{2}$. By Theorem 2.2 (1), one has
$$E_{p^{k}-1}(1,x)=E_{p^{k}-1}(1,y(1-y))=\frac{y^{p^{k}}-(1-y)^{p^{k}}}{2y-1}
=\frac{(2y-1)^{p^{k}}}{2y-1}$$
$$=(2y-1)^{p^{k}-1}=[(2y-1)^{2}]^{\frac{p^{k}-1}{2}}
=[-4y(1-y)+1]^{\frac{p^{k}-1}{2}}=(-4x+1)^{\frac{p^{k}-1}{2}}.$$
Also Theorem 2.2 (1) implies that
$$E_{p^{k}-1}\Big(1,\frac{1}{4}\Big)=\frac{p^k}{2^{p^k-1}}=0
=\Big(1-4\times\frac{1}{4}\Big)^{\frac{p^{k}-1}{2}}$$
as one desires.
\end{proof}

\noindent{\bf Lemma 2.1.} \cite{[LN]} {\it Each of the following is true:\\
{\rm(1).} Every linear polynomial over $\mathbb{F}_{q}$ is a PP of $\mathbb{F}_{q}$.\\
{\rm(2).} The monomial $x^n$ is a PP of $\mathbb{F}_{q}$ if and only if $(n, q-1)=1$.}

By Proposition 2.3 and Lemma 2.1, the following result follows immediately.

\noindent{\bf Corollary 2.1.} {\it Let $p$ be an odd
prime and $q=p^e$. Let $e$ and $k$ be positive integers with $1\leq k\leq e$. Then
$E_{p^{k}-1}(1,x)$ is a PP of $\mathbb{F}_{q}$ if and only if
$(\frac{p^{k}-1}{2},q-1)=1.$}

\noindent{\bf Lemma 2.2.} \cite{[HMSY]} {\it Let $x\in\mathbb{F}_{q^2}$. Then
$x(1-x)\in\mathbb{F}_q$ if and only if $x^q=x$ or $x^q=1-x$.}

We define
$$V:=\{x\in \mathbb{F}_{q^2}: x^q=1-x\}.$$
Then $\mathbb{F}_q\cap V=\{\frac{1}{2}\}$. We can now give a
characterization for $E_{n}(1,x)$ to be a PP.

\noindent{\bf Theorem 2.3.} {\it Let $p$ be an odd prime and
$f: y\mapsto \frac{y^{n+1}-(1-y)^{n+1}}{2y-1}$
be a mapping on $(\mathbb{F}_q\cup V)\setminus \{\frac{1}{2}\}$.
Then $E_{n}(1,x)$ is a PP of $\mathbb{F}_q$ if and only
if $f$ is $2$-to-$1$ and $f(y)\neq \frac{n+1}{2^{n}}$ for any $y\in
(\mathbb{F}_q\cup V)\setminus \{\frac{1}{2}\}$.}
\begin{proof}
First to show the sufficiency part, we choose two elements
$x_{1}$ and $x_{2}\in\mathbb{F}_q$ satisfying that
$E_{n}(1,x_{1})=E_{n}(1,x_{2})$. Since $x_1, x_2\in\mathbb{F}_q$, there exist
$y_1, y_2\in\mathbb{F}_{q^2}$ such that
$x_{1}=y_{1}(1-y_{1})$ and $x_{2}=y_{2}(1-y_{2})$.
Then by Lemma 2.2, we know that $y_{1}, y_2\in \mathbb{F}_q\cup V$.
Consider the following cases.

{\sc Case 1.} Exactly one of $x_1$ and $x_2$ is equal to $\frac{1}{4}$.
Without loss of any generality, one may let $x_{1}=\frac{1}{4}$.
Then $y_{1}=\frac{1}{2}$. Since $E_{n}(1,x_{1})=E_{n}(1,x_{2})$,
it follows from Theorem 2.2 (1) that
$E_{n}(1,x_{2})=E_{n}(1,\frac{1}{4})=\frac{n+1}{2^{n}}$.
Claim that $x_2=\frac{1}{4}$. Otherwise, we have
$x_{2}\neq\frac{1}{4}$. It follows that $y_{2}\neq\frac{1}{2}$. Since
$f(y)\neq \frac{n+1}{2^{n}}$ for any $y\in (\mathbb{F}_q\cup V)\setminus \{\frac{1}{2}\}$,
by Theorem 2.2 (1) we derive that $$E_{n}(1,x_{2})=E_{n}(1,y_{2}(1-y_{2})
=\frac{y_{2}^{n+1}-(1-y_{2})^{n+1}}{2y_{2}-1}=f(y_2)\ne \frac{n+1}{2^{n}},$$
which arrives at a contradiction. Hence we must have $x_{2}=\frac{1}{4}$.
The claim is proved. Now by the claim, one has $x_1=x_2$.

{\sc Case 2.} $x_{1}\neq\frac{1}{4}$ and $x_{2}\neq\frac{1}{4}$.
Since $E_{n}(1,x_{1})=E_{n}(1,x_{2})$, we have
$f(y_1)=f(y_2)$. Since $f$ is a $2$-to-$1$ mapping
on $(\mathbb{F}_q\cup V)\setminus \{\frac{1}{2}\}$,
it follows that $y_{1}=y_{2}$ or $y_{1}=1-y_{2}$. This implies that
$x_{1}=x_{2}$. Hence $E_{n}(1,x)$ is a PP of $\mathbb{F}_q$.
Therefore the sufficiency part is proved.

Let us now prove the necessity part. Assume that $E_{n}(1,x)$ is a PP
of $\mathbb{F}_q$. We choose two elements $y_{1},y_{2}\in(\mathbb{F}_q\cup  V)
\setminus \{\frac{1}{2}\}$ such that $f(y_1)=f(y_2)$, namely,
$$\frac{y_{1}^{n+1}-(1-y_{1})^{n+1}}{2y_{1}-1}=
\frac{y_{2}^{n+1}-(1-y_{2})^{n+1}}{2y_{2}-1}.\eqno(2.1)$$
Since $y_{1},y_{2}\in (\mathbb{F}_q\cup V)\setminus \{\frac{1}{2}\}$,
by Lemma 2.2 one has $y_{1}(1-y_{1})\in\mathbb{F}_q$ and
$y_{2}(1-y_{2})\in\mathbb{F}_q$. Then by Theorem 2.2 (1), (2.1) infers that
$$E_{n}(1,y_{1}(1-y_{1}))=E_{n}(1,y_{2}(1-y_{2})).$$
But $E_{n}(1,x)$ is a PP of $\mathbb{F}_q$, we then have $y_{1}(1-y_{1})=y_{2}(1-y_{2})$.
Thus one can immediately get that $y_{1}=y_{2}$ or $y_{1}=1-y_{2}$. Thus $f$ is
a 2-to-1 mapping on $(\mathbb{F}_q\cup V)\setminus \{\frac{1}{2}\}$.

Finally, picking $y\in(\mathbb{F}_q\cup V)\setminus \{\frac{1}{2}\}$,
it follows from Lemma 2.2 that $y(1-y)\in\mathbb{F}_q$ and
$y(1-y)\neq\frac{1}{2}(1-\frac{1}{2})$.
Since $E_{n}(1,x)$ is a PP of $\mathbb{F}_q$, it follows that
$$E_{n}(1,y(1-y))\neq E_{n}\Big(1,\frac{1}{2}\Big(1-\frac{1}{2}\Big)\Big).$$
Note that $E_{n}(1,\frac{1}{2}(1-\frac{1}{2}))=\frac{n+1}{2^{n}}.$
Then by Theorem 2.2 (1) one has
$$\frac{y^{n+1}-(1-y)^{n+1}}{2y-1}
\neq \frac{n+1}{2^{n}}.$$
Thus $f(y)\ne \frac{n+1}{2^{n}}$ for any $y\in
(\mathbb{F}_q\cup V)\setminus \{\frac{1}{2}\}$.
The necessity part is proved.

This completes the proof of Theorem 2.3.
\end{proof}

\section{\bf Necessary conditions for $E_{n}(1,x)$ to be permutational}

In the present section, we study some necessary conditions on $n$ for
$E_{n}(1,x)$ to be a PP of $\mathbb{F}_q$.
Note that $E_{n}(1,0)=1$. By the following recursive relation
\[\begin{cases} E_{0}(1,1)=1,\\
E_{1}(1,1)=1,\\
E_{n+2}(1,1)=E_{n+1}(1,1)-E_{n}(1,1),
\end{cases}\]
it follows that
$$E_{2}(1,1)=0, E_{3}(1,1)=-1, E_{4}(1,1)=-1, E_{5}(1,1)=0.$$
The sequence $\{E_{n}(1,1)\mid n\in\mathbb{N} \}$ has period $6$ and
\[E_{n}(1,1)=
\begin{cases}\ \ 0,\ \ {\rm if} \ n\equiv 2,5\pmod 6;\\
\ \ 1,\ \ {\rm if} \ n\equiv 0,1\pmod 6;\\
-1,\ \ {\rm if} \ n\equiv 3,4\pmod 6.
\end{cases}\]
\noindent{\bf Theorem 3.1.} {\it Assume that $E_{n}(1,x)$ is a PP of
$\mathbb{F}_q$. If $p=2$, then $3\mid (n+1)$; If $p$ is an odd
prime, then $n\not\equiv 0,1(\textup{mod}\ 6)$.}
\begin{proof}
By comparing $E_{n}(1,0)$ with $E_{n}(1,1)$, we get the desired
result immediately.
\end{proof}

Let $m\geq0$ be an integer. We define the polynomial $f_m(x)$ by
$$f_{m}(x):=\sum_{j=0}^{[\frac{m-1}{2}]}\binom{m}{2j+1}x^{j}\in \mathbb{Z}[x].$$
We have the following relation between $f_{n+1}(x)$ and $E_{n}(1,x)$.

\noindent{\bf Theorem 3.2.} {\it Let $p$ be an odd prime. Then
$E_{n}(1,x)=\frac{1}{2^{n}}f_{n+1}(1-4x).$ Consequently,
$E_{n}(1,x)$ is a PP of $\mathbb{F}_q$ if and only if $f_{n+1}(x)$
is a PP of $\mathbb{F}_q$.}

\begin{proof}
First we write $x=y(1-y)$ with an indeterminate $y\ne \frac{1}{2}$.
Let $u=2y-1$. Then by Theorem 2.2 (1), we derive that

\begin{eqnarray*}
 &&E_{n}(1,x)=E_{n}(1,y(1-y))\\
&=&\frac{1}{u}[y^{n+1}-(1-y)^{n+1}]\\
&=&\frac{1}{u}\Big[\Big(\frac{1+u}{2}\Big)^{n+1}-\Big(\frac{1-u}{2}\Big)^{n+1}\Big]\\
&=&\frac{1}{2^{n+1}u}\big[(1+u)^{n+1}-(1-u)^{n+1}\big]\\
&=&\frac{1}{2^{n}u}\sum_{j=0}^{[\frac{n}{2}]}\binom{n+1}{2j+1}u^{2j+1}\\
&=&\frac{1}{2^{n}}\sum_{j=0}^{[\frac{n}{2}]}\binom{n+1}{2j+1}u^{2j}\\
&=&\frac{1}{2^{n}}f_{n+1}(u^{2})\\
&=&\frac{1}{2^{n}}f_{n+1}(1-4y(1-y))\\
&=&\frac{1}{2^{n}}f_{n+1}(1-4x).
\end{eqnarray*}

Next let $x=\frac{1}{4}$. Then we obtain that
$$E_{n}(1,x)=E_{n}\Big(1,\frac{1}{4}\Big)=\frac{n+1}{2^{n}}
=\frac{1}{2^{n}}f_{n+1}(0)=\frac{1}{2^{n}}f_{n+1}(1-4x).$$
So the first part is proved.

Since $\frac{1}{2^{n}}\in\mathbb{F}_{q}^*$ and $1-4x$ is linear, we know that
$E_{n}(1,x)$ is a PP of $\mathbb{F}_q$ if and only if $f_{n+1}(x)$
is a PP of $\mathbb{F}_q$. The proof of Theorem 3.2 is complete.
\end{proof}

Using the the relation between $f_{n+1}(x)$ and $E_{n}(1,x)$ is described in Theorem 3.2,
we can get the following results.

\noindent{\bf Theorem 3.3.} {\it Let $p$ be an odd prime and $m$ be a nonnegative integer
with $p\not|(m+1)$. If $E_{2m+1}(1,x)$ is a PP of $\mathbb{F}_q$, then $m$ is odd and $(m,q-1)=1$.}
\begin{proof}
We suppose that $E_{2m+1}(1,x)$ is a PP of $\mathbb{F}_q$. Then it follows from
Theorem 3.2 that $f_{2m+2}(x)$ is a PP of $\mathbb{F}_q$. So we can choose an
element $x_{0}\in\mathbb{F}_q$ such that $f_{2m+2}(x_{0})=0$. Since $f_{2m+2}(0)=2m+2\neq0$ and
$f_{2m+2}(x)$ is a PP of $\mathbb{F}_q$, we deduce that $x_{0}\neq0$.

On the other hand, one can easily check that $f_{2m+2}(x)=x^mf_{2m+2}(x^{-1})$. Namely,
$f_{2m+2}(x)$ is a self-reciprocal polynomial. Then by $f_{2m+2}(x_{0})=0$ and $x_{0}\ne0$,
we have that $f_{2m+2}(x_{0})=f_{2m+2}(x_{0}^{-1})=0$. Since $f_{2m+2}(x)$ is a PP of
$\mathbb{F}_q$, we derive that $x_{0}=x_{0}^{-1}$, i.e., $x_{0}=\pm 1$. But
$$f_{2m+2}(1)=\sum_{j=0}^m\binom{2m+2}{2j+1}=2^{2m+1}\ne 0.$$
Then $x_{0}$ must equal $-1$. Thus we have
\begin{eqnarray*}
0=f_{2m+2}(-1)&=&\sum_{j\equiv1\pmod 4}\binom{2m+2}{j}
-\sum_{j\equiv3\pmod 4}\binom{2m+2}{j}\\
 &=&\frac{1}{2}\big[i(1-i)^{2m+2}-i(1+i)^{2m+2}\big]\\
&=&\frac{1}{2}i\big[(\sqrt{2}e^{\frac{-\pi i}{4}})^{2m+2}
-(\sqrt{2}e^{\frac{\pi i}{4}})^{2m+2}\big]\\
&=&2^{m}i\big[e^{\frac{-(m+1)\pi i}{2}}-e^{\frac{(m+1)\pi
i}{2}}\big].
\end{eqnarray*}
It follows that $e^{\frac{-(m+1)\pi i}{2}}-e^{\frac{(m+1)\pi
i}{2}}=0$. Hence $m+1$ is even. In other words, $m$ is odd.

Let us show that $(m,q-1)=1$. Assume that
$(m,q-1)=d\geq3$. Let $\theta\in \mathbb{F}_{q}^{\ast}$ satisfy
$o(\theta)=d$, where $o(\theta)$ means the order of $\theta$ in
$\mathbb{F}_{q}^{\ast}$. Since $f_{2m+2}(x)$ is self-reciprocal, one has
$f_{2m+2}(\theta)=\theta^m f_{2m+2}(\theta^{-1})=f_{2m+2}(\theta^{-1})$. But
$\theta\neq\theta^{-1}$, which contradicts with the fact that $f_{2m+2}(x)$ is a
PP of $\mathbb{F}_q$. Thus $(m,q-1)=1$ as required.

This completes the proof of Theorem 3.3.
\end{proof}
The following lemmas are needed in the reminder of this section.

\noindent{\bf Lemma 3.1.} {\it Let $p$ be an odd prime and $q$ be
the power of $p$. Let $n\ge 1$ be an integer
with $n\equiv1\pmod 4$. Then $(n+1, q-1)(n+1, q+1)=2(n+1, q^2-1)$.}
\begin{proof}
Since $q$ is odd and $n\equiv1\pmod 4$, we have $(n+1, q-1, q+1)=2$.
Let $(n+1, q-1)=2d_1$ and $(n+1, q+1)=2d_2$. Then $d_1$ and $d_2$ are two odd integer,
$(d_1, d_2)=1$ and $n+1=2d_1d_2l$ for some positive integer $l$. Since $n\equiv1\pmod 4$,
it follows that $n+1\equiv2\pmod 4$ and $(l, 2)=1$. Let $q-1=2d_1u_1$ and $q+1=2d_2u_2$.
Then one can deduce that $(d_2l, u_1)=1$ and $(d_1l, u_2)=1$.
It implies that $(l, u_1)=(l, u_2)=1$. Thus $(l, 2u_1u_2)=1$. It then follows that
$$(n+1, q-1)(n+1, q+1)=4d_1d_2=4d_1d_2(l, 2u_1u_2)$$
$$=2(2d_1d_2l, 4d_1d_2u_1u_2)=2(n+1, q^2-1)$$
as desired. Lemma 3.1 is proved.
\end{proof}
\noindent{\bf Lemma 3.2. \cite{[HL]}} {\it Let $\theta\not\in \{0,1\}$ be in
some extension of $\mathbb{F}_{q}$ and let
$y=\frac{\theta+1}{\theta-1}$. Then $y^{2}\in \mathbb{F}_{q}$ if and
only if $\theta^{q+1}=1$ or $\theta^{q-1}=1$.}

\noindent{\bf Theorem 3.4.} {\it Let $p>3$ be an odd prime and $n\geq0$
be an integer with $3\mid(n+1)$ and $n\equiv1\pmod 4$. If $E_{n}(1,x)$ is a PP of
$\mathbb{F}_q$, then $(n+1,q^2-1)=6$. }
\begin{proof}
Since $p>3$, we get that $q\equiv1\ {\rm or}\ 2\pmod3$. Thus $3\mid (q+1)$ or
$3\mid (q-1)$. Namely, $3$ divides $q^2-1$.  Since $n+1$ is
divisible by $3$, we get that $3\mid (n+1,q^2-1)$. But $p$ and $n$ are
odd integers, we deduce that $2\mid (n+1,q^2-1)$. Thus $6\mid (n+1,q^2-1)$.
That is, $(n+1,q^2-1)\geq6$. In what follows we show that $(n+1,q^2-1)=6$.

Assume that $(n+1,q^2-1)>6$. Writing
$$E:=\{\theta\in \mathbb{F}_{q^2}^{\ast}: \theta\neq1,
\theta^{(n+1,q+1)}=1\ {\rm or}\ \theta^{(n+1,q-1)}=1\}$$
gives us that 
$$|E|=(n+1,q+1)+(n+1,q-1)-3.$$
Then it follows from Lemma 3.1 and the assumption $(n+1,q^2-1)>6$ that
$$(n+1,q-1)(n+1,q+1)=2(n+1,q^2-1)>12.$$
From this inequality one can derive that $|E|>4$.

We take three distinct elements $\theta_{1},\theta_{2},\theta_{3} \in E$.
Let $i$ be an integer with $1\le i\le 3$.
Then $\theta_i^{q+1}=1$ or $\theta_i^{q-1}=1$.
Let $y_{i}=\frac{\theta_{i}+1}{\theta_{i}-1}$.
It follows from Lemma 3.2 that $y_{i}^{2}\in \mathbb{F}_q$. Since
$y_{i}=\frac{\theta_{i}+1}{\theta_{i}-1}$, we have
$\frac{y_{i}+1}{y_{i}-1}=\theta_{i}$. Thus
$(\frac{y_{i}+1}{y_{i}-1})^{n+1}=1$. Namely,
$(y_{i}+1)^{n+1}=(y_{i}-1)^{n+1}$. So by
$$f_{n+1}(y_{i}^{2})=\frac{1}{2y_{i}}[(1+y_{i})^{n+1}-(1-y_{i})^{n+1}],$$
we deduce that $f_{n+1}(y_{i}^{2})=0$. Since $\theta_{1},\theta_{2},
\theta_{3} \in E$ are distinct, it is easy to check that
$y_{1},y_{2}$ and $y_{3}$ are distinct. Thus at least two of
$y_{1}^{2}, y_{2}^{2}$ and $y_{3}^{2}$ are distinct.
But
$$f_{n+1}(y_{1}^{2})=f_{n+1}(y_{2}^{2})=f_{n+1}(y_{3}^{2})=0.$$
Hence $f_{n+1}(x)$ is not a PP of
$\mathbb{F}_q$. By Theorem 3.2, one derives that $E_{n}(1,x)$ is not a
PP of $\mathbb{F}_q$. This is a contradiction. Thus $(n+1,q^2-1)=6$
as desired.

The proof of Theorem 3.4 is complete.
\end{proof}
\noindent{\bf Theorem 3.5.} {\it Let $p>3$ be an odd prime and $n\geq0$
be an integer with $3\not|(n+1)$ and $n\equiv1\pmod 4$.
If $E_{n}(1,x)$ is a PP of $\mathbb{F}_q$, then $(n+1,q^2-1)=2$.}
\begin{proof}
Since $3$ does not divide $n+1$, we have $2|(n+1,q^2-1)$.
Let us show that $(n+1,q^2-1)=2$. Assume that $(n+1,q^2-1)>2$.
Then $(n+1,q^2-1)\ge 6$. Let
$$E:=\{\theta\in \mathbb{F}_{q^2}^{\ast}: \theta\neq1,
\theta^{(n+1,q+1)}=1\ {\rm or}\ \theta^{(n+1,q-1)}=1\}.$$
Then 
$$|E|=(n+1,q+1)+(n+1,q-1)-3.$$ 
By Lemma 3.1, one has
$$(n+1,q-1)(n+1,q+1)=2(n+1,q^2-1)\ge 12.$$
Then one derives that $|E|\geq4$.
Then in the similar way as in the proof of Theorem 3.4, we can show that
$f_{n+1}(x)$ is not a PP of $\mathbb{F}_q$. Then by Theorem 3.2 we obtain
that $E_{n}(1,x)$ is not a PP of $\mathbb{F}_q$, which is a contradiction.
We then conclude that $(n+1,q^2-1)=2$.

This ends the proof of Theorem 3.5.
\end{proof}

\section{\bf Computation of $\sum_{a\in \mathbb{F}_{q}}E_{n}(1,a)$ }

In this section, we compute the sum $\sum_{a\in \mathbb{F}_{q}}E_{n}(1,a)$.
By Proposition 2.2, we have
\begin{eqnarray*}
\sum_{n\geq0}E_{n}(1,x)t^{n}&=&\frac{1}{1-t+xt^{2}}\\
&=&\frac{1}{1-t}\ \frac{1}{1-\frac{t^2x}{t-1}}\\
&=&\frac{1}{1-t}\sum_{k\geq 0}\Big(\frac{t^2}{t-1}\Big)^kx^k\\
&=&\frac{1}{1-t}\big[1+\sum_{k=1}^{q-1}\sum_{l\geq
0}\Big(\frac{t^2}{t-1}\Big)^{k+l(q-1)}x^{k+l(q-1)}\big]\\
&\equiv& \frac{1}{1-t}\big[1+\sum_{k=1}^{q-1}\sum_{l\geq
0}\Big(\frac{t^2}{t-1}\Big)^{k+l(q-1)}x^{k}\big]\pmod {x^q-x}\\
&=&\frac{1}{1-t}\big[1+\sum_{k=1}^{q-1}\frac{\big(
\frac{t^2}{t-1}\big)^k}{1-\big(\frac{t^2}{t-1}\big)^{q-1}}x^{k}\big]\\
&=&\frac{1}{1-t}\big[1+\sum_{k=1}^{q-1}
\frac{(t-1)^{q-1-k}t^{2k}}{(t-1)^{q-1}-t^{2(q-1)}}x^{k}\big]. \ \ \ \ \ \ \ \ \ \ \ \ \ \ \ \ \ \ \ \ \
\ \ \ \ \ \  (4.1)
\end{eqnarray*}

On the other hand, by Theorem 2.2 (3), we know that if $n_{1}\equiv
n_{2}\pmod {q^2-1}$, then $E_{n_{1}}(1,x)=E_{n_{2}}(1,x)$ for any
$x\in \mathbb{F}_{q}\setminus\{\frac{1}{4}\}$. It then follows that
\begin{eqnarray*}
\sum_{n\geq0}E_{n}(1,x)t^{n}&=&1+\sum_{n=1}^{q^2-1}\sum_{l\geq
0}E_{n+l(q^2-1)}(1,x)t^{n+l(q^2-1)}\\
&\equiv& 1+\sum_{n=1}^{q^2-1}E_{n}(1,x)\sum_{l\geq
0}t^{n+l(q^2-1)}\pmod {x^q-x}\\
&=&1+\frac{1}{1-t^{q^2-1}}\sum_{n=1}^{q^2-1}E_{n}(1,x)t^n. \ \ \ \ \ \ \ \ \ \ \ \ \ \ \
 \ \ \ \ \ \ \ \ \ \ \ \ \ \ \ \ \ \ \ \ \ \ \ \ \ \ \ \ (4.2)
\end{eqnarray*}

Now (4.1) together with (4.2) implies that
\begin{eqnarray*}
&&\sum_{n=1}^{q^2-1}E_{n}(1,x)t^n\\
&\equiv&(1-t^{q^2-1})\Big(\frac{1}{1-t}-1\Big)+
\frac{1-t^{q^2-1}}{1-t}\sum_{k=1}^{q-1}\frac{(t-1)^{q-1-k}t^{2k}}
{(t-1)^{q-1}-t^{2(q-1)}}x^{k}\pmod {x^q-x}\\
&=&\frac{t(1-t^{q^2-1})}{1-t}+h(t)\sum_{k=1}^{q-1}(t-1)^{q-1-k}t^{2k}x^k,
\ \ \ \ \ \ \ \ \ \ \ \ \ \ \
 \ \ \ \ \ \ \ \ \ \ \ \ \ \ \ \ \ \ \ \ \ \ \ \ \ \ \ \ \ \ \ (4.3)
\end{eqnarray*}
where
$$h(t):=\frac{t^{q^2-1}-1}{(t-1)^q-(t-1)t^{2(q-1)}}.$$
We need the following well-known result.

\noindent{\bf Lemma 4.1.} \cite{[LN]} {\it Let $u_0, u_1,..., u_{q-1}$ be the all
elements of $\mathbb{F}_{q}$. Then }
\[\sum_{i=0}^{q-1}u_i^k=
 \begin{cases} 0,\ \ \ \ \ \ {\it if} &0\leq k\leq q-2;\\
-1,\ \ \ \ {\it if} & k=q-1.
\end{cases}\]
Then by Lemma 4.1 and (4.3), we obtain that
\begin{eqnarray*}
&&\sum_{n=1}^{q^2-1}\Big(\sum_{a\in \mathbb{F}_{q}}E_{n}(1,a)\Big)t^n\\
&=&\sum_{n=1}^{q^2-1}E_{n}\Big(1,\frac{1}{4}\Big)t^n+\sum_{n=1}^{q^2-1}
\Big(\sum_{a\in \mathbb{F}_{q}\setminus\{\frac{1}{4}\}}E_{n}(1,a)\Big)t^n\\
&=&\sum_{n=1}^{q^2-1}E_{n}\Big(1,\frac{1}{4}\Big)t^n+\sum_{a\in
\mathbb{F}_{q}\setminus\{\frac{1}{4}\}}\sum_{n=1}^{q^2-1}E_{n}(1,a)t^n\\
&=&\sum_{n=1}^{q^2-1}\frac{n+1}{2^n}t^n+
\sum_{a\in\mathbb{F}_{q}\setminus\{\frac{1}{4}\}}\frac{t(1-t^{q^2-1})}{1-t}
+h(t)\sum_{k=1}^{q-1}(t-1)^{q-1-k}t^{2k}\sum_{a\in\mathbb{F}_{q}
\setminus\{\frac{1}{4}\}}a^k\\
&=&\sum_{n=1}^{q^2-1}\frac{n+1}{2^n}t^n+(q-1)\frac{t(1-t^{q^2-1})}{1-t}+h(t)
\sum_{k=1}^{q-1}(t-1)^{q-1-k}t^{2k}\sum_{a\in\mathbb{F}_{q}}a^k\\
&-&h(t)\sum_{k=1}^{q-1}(t-1)^{q-1-k}t^{2k}\Big(\frac{1}{4}\Big)^k\\
&=&\sum_{n=1}^{q^2-1}\frac{n+1}{2^n}t^n-\frac{t(1-t^{q^2-1})}{1-t}-h(t)t^{2(q-1)}-h(t)
\sum_{k=1}^{q-1}(t-1)^{q-1-k}t^{2k}\Big(\frac{1}{4}\Big)^k. \ \ \ \ \ (4.4)
\end{eqnarray*}

However, we have
$$h(t)=\frac{t^{q^2-1}-1}{(1-t^{q-1})(t^{q}-t^{q-1}-1)}
=\frac{t^{q^2}-t}{(t-t^{q})(t^{q}-t^{q-1}-1)}
:=\frac{\sum_{i=0}^{q^2-q}b_it^i}{t^{q}-t^{q-1}-1}.\eqno(4.5)$$
Evidently, $\sum_{i=0}^{q^2-q}b_it^i=-1-(t-t^q)^{q-1}$. Then
the binomial expansion theorem applied to $(t-t^q)^{q-1}$
gives us the following result.

\noindent{\bf Proposition 4.1. } {\it For $0\leq i\leq q^2-q$,
write $i=\alpha+\beta q$ with $0\leq \alpha, \beta\leq q-1$. Then}
\[b_i=
\begin{cases} (-1)^{\beta+1}\binom{q-1}{\beta},\ \ \ {\it if}\ \alpha+\beta=q-1;\\
 -1,\ \ \ \ \ \ \ \ \ \ \ \ \ \ \ \ \ {\it if}\ \alpha=\beta=0;\\
 0, \ \ \ \ \ \ \ \ \ \ \ \ \ \ \ \ \ \ \ \  {\it otherwise}.\\
\end{cases}\]

Let $a_n:=\sum_{a\in\mathbb{F}_{q}}E_{n}(1,a)$ for $1\le n\le q^2-1$.
Then by (4.4) and (4.5), we arrive at
$$\sum_{n=1}^{q^2-1}\Big(a_n-\frac{n+1}{2^n}\Big)t^n=-\frac{t(1-t^{q^2-1})}{1-t}
-\frac{\sum_{i=0}^{q^2-q}b_it^i}{t^{q}-t^{q-1}-1}\Big(t^{2(q-1)}+
\sum_{k=1}^{q-1}(t-1)^{q-1-k}t^{2k}\Big(\frac{1}{4}\Big)^k\Big).$$
It infers that

\begin{eqnarray*}
&&(t^{q}-t^{q-1}-1)\sum_{n=1}^{q^2-1}\Big(a_n-\frac{n+1}{2^n}\Big)t^n\\
&=&(1-t^{q}+t^{q-1})\sum_{i=1}^{q^2-1}t^i-\Big(t^{2(q-1)}+\sum_{k=1}
^{q-1}(t-1)^{q-1-k}t^{2k}\big(\frac{1}{4}\big)^k\Big)
\Big(\sum_{i=0}^{q^2-q}b_it^i\Big). \ \ \ \ \ \ \ \ (4.6)
\end{eqnarray*}

We let $\sum_{i=1}^{q^2+q-1}c_it^i$ denote the right-hand side of (4.6)
and write $d_n:=a_n-\frac{n+1}{2^n}$ for integer $n$ with $1\le n\le q^2-1$.
Then (4.6) tells us that
$$(t^{q}-t^{q-1}-1)\sum_{n=1}^{q^2-1}d_nt^n=\sum_{i=1}^{q^2+q-1}c_it^i.\eqno(4.7)$$
By comparing the coefficient of $t^i$ with $1\le i \le q^2+q-1$ in both sides of (4.7),
one obtains the following recursive relations:
\[\begin{cases} c_j=-d_j,\ \ \ \ \ \ \ \ \ \ \ \  \ \ \ \ \ \ \  \ \
\ \ \ \ \  \ \ \ \ \ \ {\rm if} \ 1\leq j\leq q-1;\\
c_{q}=-d_{1}-d_{q}; \ \ \ \ \ \ \ \  \ \ \ \ \ \ \ \ \ \  \\
c_{q+j}=d_{j}-d_{j+1}-d_{q+j}, \ \ \ \ \  \ \ \ \ \ \ \ {\rm if} \ 1\leq j\leq q^2-q-1;\\
c_{q^2+j}=d_{q^2-q+j}-d_{q^2-q+j+1},  \ \ \  \ \ \ {\rm if} \ 0\leq j\leq q-2;\\
c_{q^2+q-1}=d_{q^2-1}.\\
\end{cases}\]
It then follows that
\[\begin{cases} d_j=-c_j,\ \ \ \ \ \ \ \ \ \ \ \  \ \ \ \ \ \ \  \ \ \ \
\ \ \  \ \ \ \ \ \ \ \  \ \ \ \ \ \ \  \  {\rm if} \ 1\leq j\leq q-1;\\
d_{q}=c_{1}-c_{q}; \ \ \ \ \ \ \ \  \ \ \ \ \ \ \ \ \ \  \\
d_{lq+j}=d_{(l-1)q+j}-d_{(l-1)q+j+1}-c_{lq+j},\ \ \
{\rm if}\ 1\leq l \leq q-2\ {\rm and} \ 1\leq j\leq q-1;\\
d_{lq}=d_{(l-1)q}-d_{(l-1)q+1}-c_{lq},\ \ \ \ \ \ \ \ \ \ \ \ \ \ \ {\rm if}\ 2\leq l \leq q-2;\\
d_{q^2-q+j}=\sum_{i=j}^{q-1}c_{q^2+i},  \ \ \ \ \ \ \ \ \ \ \ \ \ \ \  \ \  \ \ \ \ \ \ \  \ {\rm if} \ 0\leq j\leq q-1.\\
\end{cases} \eqno(4.8)\]

One can now give the main result of this section as the conclusion of this paper.

\noindent{\bf Theorem 4.1.} {\it Let $c_i$ be given as above for
$1\le i \le q^2+q-1$. Then each of the following is true: }
\begin{eqnarray*}
&&\sum_{a\in\mathbb{F}_{q}}E_{j}(1,a)=-c_j+\frac{j+1}{2^j} \ {\it if} \ 1\leq j\leq q-1;\\
&&\sum_{a\in\mathbb{F}_{q}}E_{q}(1,a)=c_{1}-c_{q}+\frac{1}{2^q};\\
&&\sum_{a\in\mathbb{F}_{q}}E_{lq+j}(1,a)=\sum_{a\in\mathbb{F}_{q}}E_{(l-1)q+j}(1,a)
-\sum_{a\in\mathbb{F}_{q}}E_{(l-1)q+j+1}(1,a)-c_{lq+j}
-\frac{2^{q-1}j-j-1}{2^{lq+j}}\\ &\ &\ {\it if} \ 1\leq l \leq q-2 \ {\it and} \ 1\leq j\leq q-1;\\
&&\sum_{a\in\mathbb{F}_{q}}E_{lq}(1,a)=\sum_{a\in\mathbb{F}_{q}}E_{(l-1)q}(1,a)
-\sum_{a\in\mathbb{F}_{q}}E_{(l-1)q+1}(1,a)-c_{lq}+\frac{1}{2^{lq}}\ \ {\it if} \ 2\leq l \leq q-2;\\
&&\sum_{a\in\mathbb{F}_{q}}E_{q^2-q+j}(1,a)=\sum_{i=j}^{q-1}c_{q^2+i}
+\frac{j+1}{2^{q^2-q+j}} \ {\it if} \ 0\leq j\leq q-1.
\end{eqnarray*}
\begin{proof}
Since $\sum_{a\in\mathbb{F}_{q}}E_{n}(1,a)=d_n+\frac{n+1}{2^n}$,
then by (4.8), Theorem 4.1 follows immediately.
\end{proof}

\end{document}